\documentclass{article}
\usepackage{amsmath}
\usepackage{mathrsfs}
\usepackage{txfonts}
\usepackage{bm}
\usepackage{latexsym,amsfonts,amssymb}
\begin{document}
\setcounter{page}{1}
\newtheorem{thm}{Theorem}[section]
\newtheorem{fthm}[thm]{Fundamental Theorem}
\newtheorem{dfn}[thm]{Definition}
\newtheorem{rem}[thm]{Remark}
\newtheorem{lem}[thm]{Lemma}
\newtheorem{cor}[thm]{Corollary}
\newtheorem{exa}[thm]{Example}
\newtheorem{pro}[thm]{Proposition}
\newtheorem{prob}[thm]{Problem}

\newtheorem{con}[thm]{Conjecture}
\renewcommand{\thefootnote}{\fnsymbol{footnote}}
\newcommand{\qed}{\hfill\Box\medskip}
\newcommand{\proof}{{\it Proof.\quad}}
\newtheorem{ob}[thm]{Observation}
\newcommand{\rmnum}[1]{\romannumeral #1}
\renewcommand{\abovewithdelims}[2]{%
\genfrac{[}{]}{0pt}{}{#1}{#2}}

\renewcommand{\thefootnote}{\fnsymbol{footnote}}

\title{\bf  The endomorphisms  of Grassmann graphs\footnote{Supported
by National Natural Science Foundation of China (Projects 11371072, 11301270, 11271047, 11371204).}}

\author{Li-Ping Huang\textsuperscript{a}\quad Benjian Lv\textsuperscript{b}\footnote{Corresponding author.\newline
E-mail address: lipingmath@163.com (L. Huang), benjian@mail.bnu.edu.cn (B. Lv),  wangks@bnu.edu.cn (K. Wang)}\quad Kaishun
Wang\textsuperscript{b}\\
{\footnotesize  \textsuperscript{a} \em  School of Mathematics,
Changsha University of Science and
Technology, Changsha, 410004, China}\\
{\footnotesize  \textsuperscript{b} \em Sch. Math. Sci. {\rm \&}
Lab. Math. Com. Sys., Beijing Normal University, Beijing, 100875,
China } }
 \date{}
 \maketitle

\begin{abstract}
 A graph $G$ is a core if  every  endomorphism of $G$  is an automorphism.
A graph  is called  a  pseudo-core   if every  its endomorphism is either an automorphism or a colouring.
Suppose that $J_q(n,m)$ is a Grassmann graph over a finite field with $q$ elements.
 We show that every  Grassmann graph  is a pseudo-core.
Moreover, $J_2(4,2)$ is not a core and $J_q(2k+1,2)$ ($k\geq 2$) is a core.
Further, if $m$ and $n-m+1$ are not relatively prime, then  $J_q(n,m)$ is a core when  $q$  is
  a sufficiently large integer.

 \vspace{2mm}


\end{abstract}

\section{Introduction}

Throughout this paper, all graphs  are finite undirected graphs
without loops or multiple edges. For a graph $G$,  we let $V(G)$ denote the vertex set of $G$. If $xy$ is an edge of $G$,  $x$ and
$y$ are said to be {\em adjacent}, denoted by $x\sim y$.  Let $G$ and $H$ be two graphs. A {\em homomorphism} $\varphi$ from $G$ to $H$ is a
mapping $\varphi: V(G)\rightarrow V(H)$ such that $\varphi(x)\sim \varphi(y)$   whenever $x\sim y$. If $H$ is the complete graph $K_r$,
then  $\varphi$ is a {\em $r$-colouring} of $G$ ({\em colouring} for short). An {\em isomorphism} from $G$ to $H$ is a bijection $\varphi: V(G)\rightarrow V(H)$
such that $x\sim y$ $\Leftrightarrow$ $\varphi(x)\sim \varphi(y)$. Graphs $G$ and $H$ are called isomorphic if there is an isomorphism from $G$ to $H$,
and denoted by $G\cong H$.
A homomorphism (resp. isomorphism) from $G$ to itself is called an {\em endomorphism} (resp. {\em automorphism}) of $G$.

Recall that a graph $G$ is a {\em core}  if  every  endomorphism of $G$  is an automorphism. A subgraph $H$ of $G$ is  a
{\em core} of $G$  if it is a core and there exists a homomorphism from $G$ to $H$. Every graph has a core, which is an induced
subgraph and is unique up to isomorphism \cite{Godsil}. A graph   is called {\em core-complete}  if it is a core or its core is complete.

A graph $G$ is called  a {\em pseudo-core}  if every  endomorphism of $G$ is either an automorphism or a colouring.
Every core is a pseudo-core. Any pseudo-core is core-complete but  not vice
versa. For more information, see \cite{Cameron,Godsil2,Liping-2014-2}.

 For a graph $G$, an important and difficult problem is to distinguish whether $G$ is a core \cite{Cameron,Godsil,Godsil2,G.Hahn,Knauer,Roberson}.
If $G$ is not a core or we  don't know  whether it is a core, then we  need to judge whether it is a pseudo-core
because the concept of pseudo-core  is the most close to the core. Recently, Godsil and Royle \cite{Godsil2}  discussed some properties of the pseudo-core
of a graph. Cameron and Kazanidis \cite{Cameron}  discussed the core-complete graph and the cores of symmetric graphs.
The literature \cite{Huang-Li} showed that every bilinear forms graph is a pseudo-core which is not a core.
One of the latest result is that
the literature \cite{Liping-2014-2}  proved that every alternating forms graph is a pseudo-core.
Moreover, Orel \cite{Orel-1,Orel-2} proved that each symmetric bilinear forms graph (whose diameter is greater than $2$) is a core
and each Hermitian forms graph  is a  core.

Suppose that $\mathbb{F}_q$ is the finite field with $q$ elements, where $q$ is a power of a prime.
Let $V$ be an $n$-dimensional  row vector space over  $\mathbb{F}_q$ and let ${V\brack m}$ be the set of all $m$-dimensional subspaces of $V$.
 The {\em Grassmann graph}  $J_q(n,m)$ has the
vertex set ${V\brack m}$, and two
vertices are adjacent if  their intersection is  of dimension $m-1$.
If $m=1$, we have a complete graph and hence it is a core. Since $J_q(n,m)\cong J_q(n,n-m)$,  we always
assume that $4\leq 2m\leq n$  in our discussion unless specified otherwise. The number of vertices of  $J_q(n,m)$ is the Gaussian
binomial coefficient:
\begin{equation}\label{number-vertices}
{n\brack m}=\prod_{i=1}^m\frac{q^{n+1-i}-1}{q^i-1}.
\end{equation}
For $J_q(n,m)$, the distance of two vertices $X$
and $Y$ is $d(X,Y):=m-\dim (X\cap Y)$. Any Grassmann graph  is  distance-transitive \cite[Theorem 9.3.3]{Brouwera2} and  connected.
By \cite[Corollary 4.2]{Godsil2}, every distance-transitive graph is core-complete, thus every  Grassmann graph  is core-complete.
The Grassmann graph plays an important role in geometry, graph theory, association schemes and coding theory.

Recall that an  {\em  independent set} of a graph $G$ is a  set of vertices that induces an empty subgraph. The size of the largest
independent set is called the {\em independence number} of $G$, denoted by $\alpha(G)$.
The {\em chromatic number} $\chi(G)$ of $G$ is the least value of $k$ for which $G$ can be $k$-colouring.
A {\em clique} of a graph $G$ is a  complete subgraph of $G$. A clique $C$ is maximal if there is no clique of $G$ which properly contains $C$ as a subset.
A {\em maximum clique} of $G$ is a clique with the maximum size. The {\em clique number} of $G$   is the number of vertices in a  maximum clique, denoted by $\omega(G)$.

By \cite[p.273]{Godsil2}, if $G$ is a distance-transitive graph and $\chi(G)> \omega(G)$, then $G$ is a core. Unluckily,
applying the eigenvalues or the known results of graph theory for Grassmann graph,  to prove the inequality $\chi(G)> \omega(G)$ is difficult.
Thus, it is a difficult problem to verify a Grassmann graph being a core.
However, there are some Grassmann graphs which are not cores (see Section 4). Therefore, we  need to judge whether a Grassmann graph is a pseudo-core.
 So far, this is an open problem.  We solve this problem as follows:

\begin{thm}\label{mainresult}
 Every  Grassmann graph $J_q(n,m)$ is a pseudo-core.
\end{thm}

The  paper is organized as follows. In section 2, we give some properties of the maximal cliques of Grassmann graphs. In section 3, we
shall prove Theorem \ref{mainresult}. In Section 4, we  discuss cores on Grassmann graphs. We shall show that $J_2(4,2)$ is not a core,
 $J_q(2k+1,2)$ ($k\geq 2$) is a core. Moreover,  if $m$ and $n-m+1$ are not relatively prime, then  $J_q(n,m)$ is a core when  $q$  is
  a sufficiently large integer.

\section{Maximal cliques of Grassmann graph}

In this section we shall discuss some properties of the maximal cliques of Grassmann graphs.

We will denote by $|X|$  the cardinal number of a set $X$.
Suppose that $V$ is an $n$-dimensional row vector space over $\mathbb{F}_q$.
For two vector subspaces $S$ and $T$ of $V$, the {\em join} $S\vee T$  is the minimal dimensional
vector subspace containing  $S$ and $T$. We have the dimensional formula (cf. \cite[Lemma 2.1]{Huang-Diameter} or \cite{Wan}):
\begin{equation}\label{dimensionalformula}
\dim(S\vee T)=\dim(S)+\dim(T)-\dim(S\cap T).
\end{equation}

Throughout this section, suppose that  $4\leq 2m\leq n$.
For every $(m-1)$-dimensional subspace $P$
of $V$, let $[ P\rangle_m$ denote the set of all $m$-dimensional
subspaces containing $P$, which is called a {\em star}. For every
$(m+1)$-dimensional subspace $Q$ of $V$, let $\langle Q]_m$ denote
the set of all $m$-dimensional subspaces of $Q$, which is called a
{\em top}. By \cite{Chow}, every maximal clique of $J_q(n,m)$ is a
star or a top. For more information, see \cite{M.Pankov}.

By \cite[Corollary 1.9]{Wan},
\begin{equation}\label{size}
|[P\rangle_m|=\frac{q^{n-m+1}-1}{q-1},\quad | \langle
Q]_m|=\frac{q^{m+1}-1}{q-1}.
\end{equation}
If $n>2m$, every maximum clique of $J_q(n,m)$ is a star. If $n=2m$,
every maximal clique of $J_q(n,m)$ is a maximum clique. By (\ref{size}) we have
\begin{equation}\label{chromatic-number}
\mbox{$\omega(J_q(n,m))={n-m+1\brack 1}$ if $n\geq 2m$, \  or \ $\omega(J_q(n,m))={m+1\brack 1}$ if $n< 2m$. }
\end{equation}
Since $n\geq 2m$, we have
\begin{equation}\label{size00b}
\mbox{$|[P\rangle_m|\geq|\langle Q]_m|$,  \ and \ $\left|[P\rangle_m \right|> \left|\langle Q]_m \right|$ if $n>2m$.}
\end{equation}

\begin{lem}\label{Rectangular-Application1-14}
  If $[P\rangle_m\cap \langle
Q]_m\neq\emptyset,$ then the size of $[P\rangle_m\cap \langle Q]_m$
is $q+1.$
\end{lem}
\proof Since $[P\rangle_m\cap \langle Q]_m\neq\emptyset,$ one gets
$P\subseteq Q.$ It follows that $[P\rangle_m\cap \langle Q]_m$
consists of all $m$-dimensional subspaces containing $P$ in $Q$. By
\cite[Corollary 1.9]{Wan}, the desired result follows.$\qed$


\begin{lem}\label{Rectangular-Application1-13}{\rm (\cite[Corollary 4.4]{Huang-Diameter})} \
Let $\mathcal{M}_1$ and $\mathcal{M}_2$ be two distinct stars
(tops). Then $|\mathcal{M}_1\cap\mathcal{M}_2|\leq 1$.
\end{lem}

\begin{lem}\label{maximrer3422}
Suppose $[A \rangle_m\neq [B \rangle_m$.
 Then  $[A \rangle_m\cap[B \rangle_m\neq \emptyset$ if and only if   $\dim(A\cap B)=m-2$. In this case,
  $[ A \rangle_m\cap[ B \rangle_m= \{A\vee B\}$.
\end{lem}
\proof  Since $\dim(A)=\dim(B)=m-1$ and $A\neq B$, one gets  $\dim(A\vee B)\geq m$.
 If $[ A \rangle_m\cap[ B \rangle_m\neq \emptyset$, then by Lemma \ref{Rectangular-Application1-13}, there exists a vertex $C$ of $J_q(n,m)$
such that $\{C\}=[A \rangle_m\cap[B \rangle_m$. It follows from (\ref{dimensionalformula}) and $A, B\subset C$
that $C=A\vee B$ and $\dim(A\cap B)=m-2$.
Conversely, if  $\dim(A\cap B)=m-2$, then Lemma \ref{Rectangular-Application1-13} and (\ref{dimensionalformula})  imply that $C:=A\vee B$ is a vertex of $J_q(n,m)$
and hence $\{C\}=[ A \rangle_m\cap[ B \rangle_m$.
 $\qed$

\begin{lem}\label{maximreSD22N2}    Suppose $\langle P]_m\neq \langle Q]_m$.
 Then  $\langle P]_m\cap \langle Q]_m\neq \emptyset$ if and only if
 $\dim(P\cap Q)=m$. In this case,
 $\langle P]_m\cap \langle Q]_m=\{P\cap Q\}$.
\end{lem}
\proof  By $\dim(P)=\dim(Q)=m+1$ and $P\neq Q$, we have  $\dim(P\cap Q)\leq m$.
 If $\langle P]_m\cap \langle Q]_m\neq \emptyset$, then  Lemma \ref{Rectangular-Application1-13} implies that there exists a vertex $C$ of $J_q(n,m)$
such that $\{C\}=\langle P]_m\cap \langle Q]_m$.  Since  $C\subset P\cap Q$, we get
that $C=P\cap Q$ and $\dim(P\cap Q)=m$.
Conversely, if  $\dim(P\cap Q)=m$, then by $P\cap Q\in \langle P]_m\cap \langle Q]_m$ and Lemma \ref{Rectangular-Application1-13}, we have
 $\{P\cap Q\}=\langle P]_m\cap \langle Q]_m$.
 $\qed$

In the following, let $\varphi$ be an endomorphism of $J_q(n,m)$ and let ${\rm Im}(\varphi)$ be the image of $\varphi$.
\begin{lem}\label{maxclique-0a}
 If $\mathcal{M}$ is a  maximal clique, then there exists a unique maximal   clique
 containing $\varphi(\mathcal{M})$.
 \end{lem}
\proof Suppose there exist two distinct maximal cliques
$\mathcal{M}'$ and $\mathcal{M}''$ containing
$\varphi(\mathcal{M})$. Then $\varphi(\mathcal{M})\subseteq
\mathcal{M}'\cap \mathcal{M}''$.
 By Lemmas  \ref{Rectangular-Application1-14} and \ref{Rectangular-Application1-13},
  $|\mathcal{M}'\cap \mathcal{M}''|\leq q+1$.
  Since $|{\cal M}|=|\varphi({\cal
M})|$, by (\ref{size})
 we have $|\varphi(\mathcal{M})|>q+1$, a contradiction.
 $\qed$

\begin{lem}\label{maxiclique-0b} Let $\mathcal{M}$ be a  star and $\mathcal{N}$ be a top such that
$|\varphi(\mathcal{M})\cap \varphi(\mathcal{N})|> q+1. $ Then
$\varphi(\mathcal{N})\subseteq\varphi(\mathcal{M})$.
\end{lem}
\proof Let $\mathcal{N}'$ be the maximal clique containing
$\varphi(\mathcal{N})$.  Then $|\varphi(\mathcal{M})\cap
\mathcal{N}'|>q+1$. One
gets $\varphi(\mathcal{M})=\mathcal{N}'$  by
Lemmas~\ref{Rectangular-Application1-14} and
\ref{Rectangular-Application1-13}. $\qed$

\begin{lem}\label{GrassmannTheorem-02}
Suppose there exist two distinct stars  $[A\rangle_m$ and $[B\rangle_m$
such that
$$[A\rangle_m\cap [B\rangle_m=\{X\},\quad
\varphi([A\rangle_m)=\varphi([B\rangle_m).
$$
If $\varphi([A\rangle_m)$ is a star, then  $\varphi$ is a colouring of $J_q(n,m)$.
 \end{lem}
\proof Write $\mathcal{M}:=\varphi([A\rangle_m)$. Then $\varphi([B\rangle_m)=\mathcal{M}$ and $\varphi(X)\in\mathcal{M}$.
Since the restriction mapping of $\varphi$ on a maximal clique
is injective and (\ref{size00b}),
it is easy to see that $\mathcal{M}$ is a star. If ${\rm Im}(\varphi)=\mathcal{M}$, then   $\varphi$ is a colouring of $J_q(n,m)$.
Now we prove  ${\rm Im}(\varphi)=\mathcal{M}$ as follows.
Suppose that  $Y$ is any  vertex with $Y\sim X$.  Since $G:=J_q(n,m)$ is connected, it suffices to show that there exist two distinct stars
$[C\rangle_m$ and $[D\rangle_m$
 such that
 $$
 \{Y\}=[C\rangle_m\cap [D\rangle_m \quad and \quad
 \varphi([C\rangle_m)=\varphi([D\rangle_m)=\mathcal{M}.
 $$
In fact, if we can prove this point, then we can imply that $\varphi(Z)\in \mathcal{M}$ for all $Z\in V(G)$.
We prove it as follows.

Since $X\in\langle X\vee Y]_m\cap [ A \rangle_m\cap[ B \rangle_m$,
using Lemma \ref{Rectangular-Application1-13} we get $|\langle X\vee Y]_m\cap [ A \rangle_m\cap[ B \rangle_m|=1$.
By Lemma \ref{Rectangular-Application1-14} we obtain
$$
|\langle X\vee Y]_m\cap [ A \rangle_m|=|\langle X\vee Y]_m\cap [ B
\rangle_m|=q+1.
$$
 It follows that
 $$
 |\langle X\vee Y]_m\cap ([A
\rangle_m\cup [B \rangle_m)|=2q+1.
$$
Observe that
 $$
 \varphi (\langle X\vee Y]_m\cap ([ A \rangle_m\cup
[ B \rangle_m))\subseteq \varphi (\langle X\vee Y]_m)\cap \varphi ([
A \rangle_m\cup [ B \rangle_m)\subseteq \varphi (\langle X\vee
Y]_m)\cap \mathcal{M}.
$$
Since the restriction of $\varphi$ on a clique is injective, one
gets
$$
|\varphi(\langle X\vee Y]_m)\cap \mathcal{M} |\geq 2q+1>q+1.
$$
Thus, Lemma \ref{maxiclique-0b} implies that
\begin{equation}\label{inclusion}
\varphi(\langle X\vee Y]_m)\subseteq\mathcal{M}.
\end{equation}
 So $\varphi(Y)\in
\mathcal{M}$.
Write $C:=X\cap Y$. Since every vertex of $[C \rangle_m\setminus\{X\}$ is adjacent to $X$, by our claim  we have
$\varphi([C\rangle_m)= \mathcal{M}$.

Pick a vertex $Z$ such that $Z\sim Y$ and the distance from $X$ is
$2$. Write $D=Y\cap Z$. Since $Y\in[D \rangle_m\cap \langle X\vee Y]_m$, by Lemma~\ref{Rectangular-Application1-14} we have
$|[D\rangle_m\cap \langle X\vee Y]_m|= q+1$. It follows from
(\ref{inclusion}) that  $|\varphi([D\rangle_m)\cap \mathcal{M}|\geq q+1$.  Thus  Lemma \ref{Rectangular-Application1-13}  implies that
$\varphi([D \rangle_m)=\mathcal{M}$. Since  $\{Y\}=[C\rangle_m\cap [D\rangle_m$,
 $[C\rangle_m$ and $[D \rangle_m$ are the desired stars.
$\qed$

\section{Proof of Theorem \ref{mainresult}}

For the proof of Theorem \ref{mainresult},
we only need to consider the case $4\leq 2m\leq n$. We divide the proof of Theorem \ref{mainresult} into two cases: $n>2m$ and $n=2m.$

\begin{lem}\label{mainresult0b}
 If $n>2m$, then every  Grassmann graph $J_q(n,m)$ is a pseudo-core.
\end{lem}
\proof Suppose that $n>2m\geq 4$. Then by (\ref{size00b}), every maximum clique of $J_q(n,m)$ is a star. Let  $\varphi$ be an endomorphism  of $J_q(n, m)$. Then the
restriction of  $\varphi$ on any clique is injective, so $\varphi$
transfers stars to stars.

 Suppose  $\varphi$ is not a colouring.
  It suffices to show that
$\varphi$ is an automorphism. Write $G_r:=J_q(n, r)$, where $1\leq
r\leq m-1.$ By Lemma~\ref{GrassmannTheorem-02}, the images under $\varphi$ of any two distinct and intersecting stars  are distinct.
 Hence by Lemma~\ref{maximrer3422}, $\varphi$ induces an
endomorphism $\varphi_{m-1}$ of $G_{m-1}$ such that
$$
\varphi([A \rangle_m)=[\varphi_{m-1}(A) \rangle_m.
$$

Let $X$ be any vertex of $J_q(n,m)$. Then there exist two vertices $X'$ and $X''$
of $G_{m-1}$ such that  $X=X'\vee X''$. Then $[ X' \rangle_m\cap [
X'' \rangle_m=\{X\}$ and $\varphi(X)\in \varphi([ X' \rangle_m)\cap
\varphi([ X'' \rangle_m)$. Since $\varphi$ is not a colouring, by
Lemma~\ref{GrassmannTheorem-02}  $\varphi([ X' \rangle_m)$ and
$\varphi([ X'' \rangle_m)$ are two distinct stars. By
Lemma~\ref{Rectangular-Application1-13}, $[ \varphi_{m-1}(X')
\rangle_m\cap[ \varphi_{m-1}(X'') \rangle_m=\{\varphi(X)\}$. Thus
Lemma~\ref{maximrer3422} implies that
\begin{equation}\label{fact.}
\varphi(X)=\varphi_{m-1}(X')\vee\varphi_{m-1}(X'').
\end{equation}

When $m=2$, $G_1$ is a complete graph, hence it is a core. We next
show that $\varphi_{m-1}$ is not a colouring of $G_{m-1}$ for $m\geq
3.$ For any two vertices $A_1$ and $A_3$ of
$G_{m-1} $ at distance $2$, we claim that
$$
\varphi_{m-1}(A_1)\neq \varphi_{m-1}(A_3).
$$
There exists an $A_2\in V(G_{m-1})$ such that $A_1\sim A_2\sim A_3$. Write
$Y_1:=A_1\vee A_2$ and $Y_2:=A_2\vee A_3$. Then $Y_1\sim Y_2$, so
$\varphi(Y_1)\neq \varphi(Y_2)$. By (\ref{fact.}),
$$
\varphi(Y_1)=\varphi_{m-1}(A_1)\vee \varphi_{m-1}(A_2),\quad
\varphi(Y_2)=\varphi_{m-1}(A_2)\vee \varphi_{m-1}(A_3).
$$
Thus our claim is valid. Otherwise, one has $\varphi(Y_1)=\varphi(Y_2)$, a contradiction.

Pick  a star $\mathcal N$ of  $G_{m-1}$. Since the diameter of
$G_{m-1}$ is at least two, there exists a vertex $A_4\in
V(G_{m-1})\setminus\mathcal N$ that is adjacent to some vertex in $\mathcal N$. If $B\in\mathcal N$ such that $A_4$ is not  adjacent to $B$,
then $d(A_4,B)=2$. By our claim, $\varphi_{m-1}(A_4)\neq \varphi(B)$ and hence
$\varphi_{m-1}(A_4)\not\in \varphi_{m-1}(\mathcal N)$. Therefore, $\varphi_{m-1}$ is not a colouring.

By induction, we may obtain induced endomorphism $\varphi_r$ of
$G_r$ for each $r$. Furthermore,
\begin{equation}\label{fact..}
\varphi(X)=\varphi_{k_1}(X_{k_1})\vee
\varphi_{k_2}(X_{k_2})\vee\cdots\vee \varphi_{k_s}(X_{k_s}),
\end{equation}
where $X=X_{k_1}\vee X_{k_1}\vee\cdots\vee X_{k_s}\in V(G_m)$ and $1\leq \dim
(X_{k_i})=k_i\leq m-1$.

In order to show that $\varphi$ is an automorphism, it suffices to
show that $\varphi$ is injective. Assume that $X$ and $Y$ are any two distinct vertices in $G_m$ with $d(X,Y)=s$. Thus ${\rm dim}(X\cap Y)=m-s$.
If $s=1,$ then $\varphi(X)\neq \varphi(Y).$
Now suppose $s\geq 2$. There are $1$-dimensional row vectors $X_i, Y_i$, $i=1,\ldots. s$, such that $X,Y$ can be written as
$X=(X\cap Y)\vee X_1\vee \cdots \vee X_s$, $Y=(X\cap Y)\vee Y_1\vee \cdots \vee Y_s$.
Let $Z=(X\cap Y)\vee X_1\vee \cdots \vee X_{s-1}\vee Y_s\in V(G_m)$. By $X\sim Z$, $\dim(\varphi(X)\vee \varphi(Z))=m+1$.
Applying (\ref{fact..}), one has that $\varphi(X)=\varphi_{m-s}(X\cap Y)\vee \varphi_1(X_1)\vee \cdots \vee\varphi_1(X_s)$,
$\varphi(Y)=\varphi_{m-s}(X\cap Y)\vee \varphi_1(Y_1)\vee \cdots \vee \varphi_1(Y_s)$ and
$\varphi(Z)=\varphi_{m-s}(X\cap Y)\vee \varphi_1(X_1)\vee \cdots \vee \varphi_1(X_{s-1})\vee\varphi_1(Y_s)$.
 Therefore, we get $\varphi(X)\vee\varphi(Z)\subseteq \varphi(X)\vee \varphi(Y)$. It follows that $\varphi(X)\neq \varphi(Y)$.
 Otherwise, one has $\varphi(X)\vee\varphi(Z)\subseteq \varphi(X)$, a contradiction to  $\dim(\varphi(X)\vee \varphi(Z))=m+1$.
 Hence, $\varphi$ is an automorphism, as desired.

By above discussion,  $J_q(n,m)$ is a pseudo-core when $n>2m$.
$\qed$

\begin{lem}\label{mainresult0c}
 If $n=2m$, then every  Grassmann graph $J_q(n,m)$ is a pseudo-core.
\end{lem}
\proof  Suppose that $n=2m\geq 4$. For a subspace $W$ of $V$, the {\em dual subspace} $W^\perp$ of $W$ in $V$ is defined by
$$
W^\perp=\{v\in V\mid wv^{\rm t}=0,\;\forall\; w\in W\},
$$
where $v^{\rm t}$ is the transpose of $v$.

For an endomorphism $\varphi$
of $J_q(2m,m)$, define the map
$$
\varphi^\perp:V(J_q(2m,m))\longrightarrow V(J_q(2m,m)),\quad A \longmapsto \varphi(A)^\perp.
$$
Then $\varphi^\perp$ is an endomorphism of
$J_q(2m,m)$. Note that $\varphi^\perp$ is an automorphism (resp. colouring) whenever $\varphi$ is an automorphism (resp. colouring).
For any maximal clique $\mathcal{M}$ of $J_q(2m,m)$,
$\varphi(\mathcal{M})$ and $\varphi^\perp(\mathcal{M})$ are of
different types.

Next we shall show that   $J_q(2m,m)$ is
a pseudo-core.

\textbf{Case 1.} \ There exist   $[A\rangle_m$ and $\langle X]_m$  such
that $[A\rangle_m\cap\langle X]_m\neq\emptyset$ and
$\varphi([A\rangle_m)$, $\varphi(\langle X]_m)$ are of the same
type.

By Lemma~\ref{Rectangular-Application1-14}, the size of
$[A\rangle_m\cap\langle X]_m$ is $q+1$. Then
$|\varphi([A\rangle_m)\cap\varphi(\langle X]_m)|\geq q+1$. Since $\varphi([A\rangle_m)$, $\varphi(\langle X]_m)$ are of the same type, by Lemma~\ref{Rectangular-Application1-13} one gets
\begin{equation}\label{new}
\varphi([A\rangle_m)=\varphi(\langle X]_m).
\end{equation}
Note that $A\subseteq X$. Pick any  $Y\in{V\brack m+1}$ satisfying $A\subseteq Y$ and $\dim(X\cap Y)=m$. Then
$\langle Y]_m\cap[A\rangle_m\neq\emptyset$.
By Lemma~\ref{Rectangular-Application1-14}
  we have $|\varphi(\langle Y]_m)\cap\varphi([A\rangle_m)|\geq q+1.$
By Lemma~\ref{Rectangular-Application1-13} and (\ref{new}) we obtain either $\varphi(\langle Y]_m)=\varphi(\langle X]_m)$
 or $\varphi(\langle Y]_m)$ and $\varphi(\langle X]_m)$ are of different types.

  \textbf{Case 1.1.} \ There exists a $Y\in{V\brack m+1}$ such that $\varphi(\langle Y]_m)$ and $\varphi(\langle X]_m)$ are of different types.
For any  $B\in{X\cap Y\brack m-1}$, we have that $B\subseteq Y$ and $B\subseteq X$.
Since $|[B\rangle_m)\cap\langle X]_m|=|[B\rangle_m)\cap\langle Y]_m|=q+1$, we have similarly
$$
    |\varphi([B\rangle_m)\cap\varphi(\langle X]_m)|\geq q+1,\quad
    |\varphi([B\rangle_m)\cap\varphi(\langle Y]_m)|\geq q+1.
$$
 Since $\varphi(\langle Y]_m)$ and $\varphi(\langle X]_m)$ are of different types, Lemma~\ref{Rectangular-Application1-13} implies that
  $\varphi([B\rangle_m)=\varphi(\langle X]_m)$ or $\varphi([B\rangle_m)=\varphi(\langle Y]_m)$ for any $B\in{X\cap Y\brack m-1}$.

    Since the size of ${X\cap Y\brack m-1}$ is at least $3$, by above discussion, there exist two subspaces $B_1,B_2\in{X\cap Y\brack m-1}$ such that
    $\varphi([B_1\rangle_m)=\varphi([B_2\rangle_m).$ Note that $[B_1\rangle_m\cap [B_2\rangle_m\neq\emptyset$ because $X\cap Y\in B_i$ ($i=1,2$).
  If $\varphi([B_1\rangle_m)$
    is a star, then $\varphi$ is a colouring by Lemma~\ref{GrassmannTheorem-02}.  Suppose $\varphi([B_1\rangle_m)$ is a top.
    Then $\varphi^\perp([B_1\rangle_m)$ is a star. By Lemma~\ref{GrassmannTheorem-02} again, $\varphi^\perp$ is a colouring. Hence,
    $\varphi$ is also a colouring.

   \textbf{Case 1.2.} \ $\varphi(\langle Y]_m)=\varphi(\langle X]_m)$ for any $Y\in{V\brack m+1}$.
Consider a star $[C\rangle_m$ where $C$ satisfies $C\subset X$ and $\dim(C\cap A)=m-2.$ Then $(A\vee C)\subseteq X$ and $\dim(A\vee C)=m$. For any $T\in[C\rangle_m,$ since $(A\vee C)\subseteq(A\vee T)$ and $m\leq\dim(A\vee T)\leq m+1$, there exists a
   subspace $W\in{V\brack m+1}$ such that $(A\vee T)\subseteq W$ and $\dim(W\cap X)\geq m$ (because $(A\vee C)\subseteq W\cap X$).

Since $T\in \langle W]_m$, $\varphi(T)\in \varphi(\langle W]_m)$. By the condition, $\varphi(\langle W]_m)=\varphi(\langle X]_m).$ Then
$\varphi(\langle W]_m)=\varphi([A\rangle_m)$ by (\ref{new}). It follows that $\varphi(T)\in \varphi([A\rangle_m)$ for all $T\in[C\rangle_m$, and so
$\varphi([C\rangle_m)\subseteq \varphi([A\rangle_m)$. Hence,
$\varphi([C\rangle_m)= \varphi([A\rangle_m)$. Since $[C\rangle_m\cap [A\rangle_m\neq \emptyset,$
similar to the proof of Case 1.1,   $\varphi$ is a colouring.

 \textbf{Case 2.} \ For any two maximal cliques  of different types containing  common vertices, their images under  $\varphi$  are of different types.

 In this case,   $\varphi$ maps the maximal cliques of the same type to the maximal cliques of the same type.

\textbf{Case 2.1.} \ $\varphi$ maps stars to stars. In this case $\varphi$ maps tops to tops by Lemmas \ref{Rectangular-Application1-14}
and \ref{Rectangular-Application1-13}.

If there exist two distinct stars $\mathcal{M}$ and
$\mathcal{M}'$ such that $\mathcal{M}\cap\mathcal{M}'\neq \emptyset$
and $\varphi(\mathcal{M})=\varphi(\mathcal{M}')$, then $\varphi$ is a colouring by Lemma~\ref{GrassmannTheorem-02}.
Now suppose $\varphi(\mathcal{M})\neq\varphi(\mathcal{M}')$ for any two distinct stars $\mathcal{M}$ and
$\mathcal{M}'$ with $\mathcal{M}\cap\mathcal{M}'\neq \emptyset$.
By Lemma~\ref{maximrer3422}, $\varphi$ induces an endomorphism $\varphi_{m-1}$
of  $J_q(2m,m-1)$ such that $\varphi([A\rangle_m)=[\varphi_{m-1}(A)\rangle_m$. By Lemma \ref{mainresult0b},
$J_q(2m,m-1)$ is a pseudo-core. Thus, $\varphi_{m-1}$ is an automorphism or a colouring.

We claim that $\varphi_{m-1}$ is an automorphism of $J_q(2m,m-1)$.
For any $C\in{V\brack m}$ and $B\in{C\brack m-1}$, since $C\in[B\rangle_m$ and  $\varphi([B\rangle_m)=[\varphi_{m-1}(B)\rangle_m$, we have $\varphi(C)\in[\varphi_{m-1}(B)\rangle_m$. Then $\varphi_{m-1}(B)\subseteq\varphi(C),$ which implies that $\varphi_{m-1}(\langle C]_{m-1})$ is a top of $J_q(2m,m-1)$.
 If $m=2,$ our claim is valid. Now suppose $m\geq 3$ and
 $\varphi_{m-1}$ is a colouring. Then ${\rm Im}(\varphi_{m-1})$ is  a star of $J_q(2m,m-1)$.
 Note that $\varphi_{m-1}(\langle C]_{m-1})\subseteq {\rm Im}(\varphi_{m-1})$ and $ |\varphi_{m-1}(\langle C]_{m-1})|>q+1,$ contradicting to  Lemma~\ref{Rectangular-Application1-14}. Hence, our claim is valid.

\textbf{Case 2.2.} \ $\varphi$ maps stars to tops. In this case $\varphi$ maps tops to stars by Lemmas \ref{Rectangular-Application1-14}
and \ref{Rectangular-Application1-13}.

Note that $\varphi^\perp$ maps stars to stars.  By
Case 2.1,    $\varphi^\perp$ is an automorphism. Hence, $\varphi$ is an automorphism.

By above discussion, we have proved  that every  Grassmann graph $J_q(2m,m)$ is a pseudo-core.
$\qed$

By Lemmas \ref{mainresult0b} and \ref{mainresult0c}, we have proved Theorem \ref{mainresult}.

\section{Cores on Grassmann graphs}

In  this section, we shall show that $J_2(4,2)$  is not a core and $J_q(n,m)$ is a core under some conditions.

It is well-known (cf. \cite[Theorem 6.10 and Corollary 6.2]{Chartrand}) that the chromatic number of $ G $ satisfies the following inequality:
\begin{equation*}\label{chromatic-1}
 \chi( G ) \geq  {\rm max} \left\{\omega( G ), \ |V( G )|/\alpha( G ) \right\}.
\end{equation*}
By \cite[Lemma 2.7.2]{Roberson}, if  $ G $ is a  vertex-transitive
graph, then
\begin{equation}\label{d34223gh}
\chi( G )\geq \frac{|V( G )|}{\alpha( G )}\geq \omega( G ).
\end{equation}

\begin{lem}\label{tte54gdg}  Let $G$ be a  Grassmann graph. Then $G$ is a core if and only if $\chi(G)> \omega(G)$. In particular,
if $\frac{|V(G)|}{\omega(G)}$ is not an integer, then $G$ is a core.
\end{lem}
\proof By \cite[Corollary 4.2]{Godsil2}, every distance-transitive graph is core-complete, thus $G$ is core-complete.
 Then,  $\chi(G)> \omega(G)$ implies that $G$ is a core.
Conversely, if $G$ is a core, then we must have $\chi(G)> \omega(G)$.  Otherwise, there exists an endomorphism $f$  of $G$ such that $f(G)$ is a
maximum clique of $G$, a contradiction to $G$ being a core. Thus,  $G$ is a core if and only if $\chi(G)> \omega(G)$.

By \cite[p.148, Remark]{Cameron}, if the core of $G$ is complete, then $|V(G)|=\omega(G)\alpha(G)$. Assume that $\frac{|V(G)|}{\omega(G)}$ is not an integer. Then
$|V(G)|\neq\omega(G)\alpha(G)$. Therefore, the core of $G$ is not complete and hence $G$ is a core.
 $\qed$

Denote by $\mathbb{F}_q^{m\times n}$  the set of  $m\times n$  matrices over $\mathbb{F}_q$ and  $\mathbb{F}_q^n=\mathbb{F}_q^{1\times n}$.
Let $G=J_q(n,m)$ where $n> m$. If $X$ is a vertex of $G$, then $X=[\alpha_1,\ldots,\alpha_m]$ is an $m$-dimensional  subspace of
the vector space $\mathbb{F}_q^n$, where $\{\alpha_1,\ldots,\alpha_m\}$ is a basis of $X$. Thus, $X$ has a {\em matrix representation}
$\scriptsize \left(
     \begin{array}{c}
       \alpha_1 \\
       \vdots \\
       \alpha_m \\
     \end{array}
   \right)\in \mathbb{F}_q^{m\times n}$  (cf. \cite{Huang-Diameter,Wan}).
For simpleness, the matrix representation of $X\in V(G)$ is also denoted by $X$.
For matrix representations $X,Y$ of two vertices $X$ and $Y$, $X\sim Y$ if and only if
$\scriptsize{\rm rank}\left(
\begin{array}{c}
X \\
Y \\
\end{array}
\right)=m+1$. Note that if $X$ is a matrix representation then $X=PX$ (as matrix representation) for any $m\times m$ invertible matrix $P$ over $\mathbb{F}_q$.
Then, $V(G)$ has a  matrix representation
$$V(G)=\left\{X: X\in \mathbb{F}_q^{m\times n}, \ {\rm rank}(X)=m \right\}.$$

 Now, we give an example of Grassmann graph which is not a core as follows.

\begin{exa}\label{frerered}
Let $G=J_2(4,2)$. Then $G$ is not a core. Moreover, $\chi(G)=\omega(G)=7$ and $\alpha(G)=5$.
\end{exa}
\proof  Applying the matrix representation of $V(G)$, $G=J_2(4,2)$ has $35$ vertices as follows:

$\scriptsize \ A_1=\left(
\begin{array}{cccc}
1 & 0 & 0 & 0 \\
0 & 1 & 0 & 0 \\
\end{array}
\right)$, \  \
$\scriptsize A_2=\left(
\begin{array}{cccc}
1 & 0 & 1 & 0 \\
 0 & 1 & 0 & 0 \\
 \end{array}
\right)$, \
$\scriptsize A_3=\left(
\begin{array}{cccc}
1 & 0 & 0 & 1 \\
 0 & 1 & 0 & 0 \\
 \end{array}
 \right)$, \ \
$\scriptsize A_4=\left(
\begin{array}{cccc}
1 & 0 & 1 & 1 \\
 0 & 1 & 0 & 0 \\
 \end{array}
 \right)$,

$\scriptsize \ A_5=\left(
\begin{array}{cccc}
1 & 0 & 0 & 0 \\
0 & 1 & 1 & 0 \\
\end{array}
\right)$, \ \
$\scriptsize A_6=\left(
\begin{array}{cccc}
1 & 0 & 0 & 0 \\
0 & 1 & 0 & 1 \\
\end{array}
\right)$, \
$\scriptsize A_7=\left(
\begin{array}{cccc}
1 & 0 & 0 & 0 \\
0 & 1 & 1 & 1\\
\end{array}
\right)$, \ \
$\scriptsize A_8=\left(
\begin{array}{cccc}
1 & 0 & 1 & 0 \\
0 & 1 & 1 & 0 \\
\end{array}
\right)$,

$\scriptsize \ A_9=\left(
\begin{array}{cccc}
1 & 0 & 0 & 1 \\
0 & 1 & 0 & 1 \\
\end{array}
\right)$, \
$\scriptsize A_{10}=\left(
\begin{array}{cccc}
1 & 0 & 1 & 0 \\
0 & 1 & 0 & 1 \\
\end{array}
\right)$,
$\scriptsize A_{11}=\left(
\begin{array}{cccc}
1 & 0 & 0 & 1 \\
0 & 1 & 1 & 0\\
\end{array}
\right)$,
$\scriptsize A_{12}=\left(
\begin{array}{cccc}
1 & 0 & 0 & 1\\
0 & 1 & 1 & 1 \\
\end{array}
\right)$,

$\scriptsize A_{13}=\left(
\begin{array}{cccc}
1 & 0 & 1 & 0 \\
0 & 1 & 1 & 1 \\
\end{array}
\right)$,
$\scriptsize A_{14}=\left(
\begin{array}{cccc}
1 & 0 & 1 & 1 \\
0 & 1 & 0 & 1 \\
\end{array}
\right)$,
$\scriptsize A_{15}=\left(
\begin{array}{cccc}
1 & 0 & 1 & 1 \\
0 & 1 & 1 & 0\\
\end{array}
\right)$,
$\scriptsize A_{16}=\left(
\begin{array}{cccc}
1 & 0 & 1 & 1\\
0 & 1 & 1 & 1 \\
\end{array}
\right)$,


$\scriptsize A_{17}=\left(
\begin{array}{cccc}
0 & 0 & 1 & 0 \\
0 & 0 & 0 & 1 \\
\end{array}
\right)$,
$\scriptsize A_{18}=\left(
\begin{array}{cccc}
1 & 0 & 1 & 0 \\
0 & 0 & 0 & 1 \\
\end{array}
\right)$,
$\scriptsize A_{19}=\left(
\begin{array}{cccc}
0 & 1 & 1 & 0 \\
0 & 0 & 0 & 1\\
\end{array}
\right)$,
$\scriptsize A_{20}=\left(
\begin{array}{cccc}
0 & 0 & 1 & 0\\
1 & 0 & 0 & 1 \\
\end{array}
\right)$,

$\scriptsize A_{21}=\left(
\begin{array}{cccc}
0 & 0 & 1 & 0 \\
0 & 1 & 0 & 1 \\
\end{array}
\right)$,
$\scriptsize A_{22}=\left(
\begin{array}{cccc}
1 & 1 & 1 & 0 \\
0 & 0 & 0 & 1 \\
\end{array}
\right)$,
$\scriptsize A_{23}=\left(
\begin{array}{cccc}
1 & 0 & 1 & 0 \\
1 & 0 & 0 & 1\\
\end{array}
\right)$,
$\scriptsize A_{24}=\left(
\begin{array}{cccc}
0 & 0 & 1 & 0\\
1 & 1 & 0 & 1 \\
\end{array}
\right)$,

$\scriptsize A_{25}=\left(
\begin{array}{cccc}
0 & 1 & 1 & 0 \\
0 & 1 & 0 & 1 \\
\end{array}
\right)$,
$\scriptsize A_{26}=\left(
\begin{array}{cccc}
1 & 1 & 1 & 0\\
1 & 1 & 0 & 1 \\
\end{array}
\right)$,
$\scriptsize A_{27}=\left(
\begin{array}{cccc}
0 & 1 & 0 & 0 \\
0 & 0 & 1 & 0 \\
\end{array}
\right)$,
$\scriptsize A_{28}=\left(
\begin{array}{cccc}
1 & 1 & 0 & 0 \\
0 & 0 & 1 & 0 \\
\end{array}
\right)$,

$\scriptsize A_{29}=\left(
\begin{array}{cccc}
0 & 1 & 0 & 0 \\
0 & 0 & 1 & 1 \\
\end{array}
\right)$,
$\scriptsize A_{30}=\left(
\begin{array}{cccc}
1 & 0 & 0 & 0 \\
0 & 0 & 0 & 1 \\
\end{array}
\right)$,
$\scriptsize A_{31}=\left(
\begin{array}{cccc}
0 & 1 & 0 & 0 \\
0 & 0 & 0 & 1 \\
\end{array}
\right)$,
$\scriptsize A_{32}=\left(
\begin{array}{cccc}
1 & 0 & 0 & 0 \\
0 & 0 & 1 & 0 \\
\end{array}
\right)$,

$\scriptsize A_{33}=\left(
\begin{array}{cccc}
1 & 1 & 0 & 0 \\
0 & 0 & 0 & 1 \\
\end{array}
\right)$,
$\scriptsize A_{34}=\left(
\begin{array}{cccc}
1 & 1 & 0 & 0 \\
0 & 0 & 1 & 1 \\
\end{array}
\right)$,
$\scriptsize A_{35}=\left(
\begin{array}{cccc}
1 & 0 & 0 & 0 \\
0 & 0 & 1 & 1 \\
\end{array}
\right)$.

Suppose that  $\mathcal{L}_1 = \{A_1, A_{10}, A_{12}, A_{15}, A_{17}\}$, $\mathcal{L}_2 = \{A_2, A_{6}, A_{20}, A_{19}, A_{34}\}$,
$\mathcal{L}_3 = \{A_3, A_{8}, A_{21}, A_{22}, A_{35}\}$, $\mathcal{L}_4 = \{A_5, A_{9}, A_{18}, A_{24}, A_{29}\}$,
$\mathcal{L}_5 = \{A_7, A_{14}, A_{23}, A_{27}, A_{33}\}$,
$\mathcal{L}_6 = \{A_4, A_{13}, A_{25}, A_{28}, A_{30}\}$, and
$\mathcal{L}_7=\{A_{11}, A_{16}, A_{26}, A_{31}, A_{32}\}$.
It is easy to see that $V(G)=\mathcal{L}_1\cup \mathcal{L}_2 \cup\cdots\cup \mathcal{L}_7$ and $\mathcal{L}_1,\ldots,\mathcal{L}_7$ are  independent sets.
Thus $\chi(G)\leq 7$. On the other hand, (\ref{d34223gh}) implies that $\chi(G)\geq \omega(G)= 7$. Therefore,  $\chi(G)= \omega(G)=7$.
It follows from Corollary \ref{tte54gdg} that $G$ is not a core. By (\ref{d34223gh}) again, we have $\alpha(G)=5$.
$\qed$

We guess that $J_q(2k,2)$ ($k\geq 2$) is not a core for all $q$ (which is a power of a prime). But this a difficult problem.
Next, we give some examples of Grassmann graph which is  a core.

\begin{exa}\label{3indepmann4}\ If   $k\geq 2$, then  $J_q(2k+1,2)$  is core.
\end{exa}
\proof  When $k\geq 2$, let  $G=J_q(2k+1,2)$. Applying (\ref{number-vertices}) and (\ref{chromatic-number}) we have
$$\frac{|V(G)|}{\omega(G)}=\frac{q^{2k+1}-1}{q^2-1}=\frac{q^{2k+1}-q}{q^2-1}+\frac{1}{q+1}.$$
Thus  $\frac{|V(G)|}{\omega(G)}$ is not an integer for any $q$ (which is a power of a prime). By Lemma \ref{tte54gdg}, $G$ is  a core.
$\qed$

Denote by $\mathbb Z$ the integer ring and $\mathbb Z[x]$ the polynomial ring in an indeterminate $x$ over $\mathbb Z$.
Let $\Phi_t(x)$ be the $t$th {\em cyclotomic polynomial} defined by
$$ \Phi_t(x)=\prod_{1\leq j\leq t \atop
                      {\rm gcd}(j, \,t)=1
                    }(x-\zeta_t^j),
$$
where $\zeta_t$ is the $t$th root of unity and ${\rm gcd}(j, \, t)$ is
the greatest common divisor of $j$ and $t$.
Recall that $\Phi_t(x)$ is an irreducible polynomial over $\mathbb Z$.
The polynomial $x^n-1$ over $\mathbb Z$ has the following factorization into irreducible polynomials over $\mathbb{Z}$:
\begin{equation}\label{factorization1}
x^n-1=\prod_{j\mid n}\Phi_j(x).
\end{equation}

In 1989, Knuth and Wilf  gave a factorization of Gaussian binomial coefficient (as a polynomials in $\mathbb{Z}[q]$) (cf. \cite{W.Y.C.Chen,knuth}):
\begin{equation}\label{Gaussianbinomialf523}
{n\brack m}=\prod_{i=1}^n(\Phi_i(q))^{\lfloor n/i\rfloor -\lfloor m/i\rfloor -\lfloor (n-m)/i\rfloor},
\end{equation}
where $\lfloor a\rfloor$ is the largest integer no more than $a$. Note that $\lfloor n/i\rfloor -\lfloor m/i\rfloor -\lfloor (n-m)/i\rfloor$ is equal to $0$ or $1$.


Write $G:=J_q(n,m)$ (where $4\leq 2m\leq n$) and $h(q):=\frac{|V(G)|}{\omega(G)}={n\brack m}/\omega(G)$, where $h(q)$ is seen as a polynomial
in an indeterminate $q$ over the rational number field. By (\ref{Gaussianbinomialf523}) one gets
that
$$
\omega(G)={n-m+1\brack 1}=\prod_{j=1}^{n-m+1}(\Phi_j(q))^{\lfloor (n-m+1)/j\rfloor -\lfloor 1/j\rfloor -\lfloor (n-m)/j\rfloor},
$$
\begin{align}\label{h(q)define}
h(q)=\prod_{j=2}^{n-m+1}(\Phi_j(q))^{\lfloor n/j\rfloor-\lfloor
m/j\rfloor-\lfloor(n-m+1)/j\rfloor}\prod_{j=n-m+2}^n(\Phi_j(q))^{\lfloor
n/j\rfloor-\lfloor m/j\rfloor-\lfloor(n-m)/j\rfloor}.
\end{align}

\begin{thm}\label{GrassmannTheorem2}
Assume that $m$ and $n-m+1$ are not relatively prime. If $q$ (which is a power of a prime) is a sufficiently large integer (i.e., there is a fixed positive
integer $c_{n,m}$ such that $q\geq c_{n,m}$), then the  Grassmann graph $J_q(n,m)$ is a core.
\end{thm}
\proof
 Note that $\lfloor x+y \rfloor-\lfloor x \rfloor -\lfloor y \rfloor$ is equal to $0$ or $1$
for all real numbers $x$ and $y$. We have that $\lfloor m/j\rfloor+\lfloor(n-m+1)/j\rfloor$ equals $\lfloor(n+1)/j\rfloor$ or $\lfloor(n+1)/j\rfloor+1$. Thus,
 $$-1\leq\lfloor n/j\rfloor-\lfloor m/j\rfloor-\lfloor(n-m+1)/j\rfloor\leq 0, \ j=2, \ldots,n-m+1.$$
Taking the greatest common factor $i$ ($i\geq2$) of $m$ and $n-m+1$. It is easy to see that
$$\lfloor n/i\rfloor-\lfloor m/i\rfloor-\lfloor(n-m+1)/i\rfloor=-1.$$
Let $f(q)=\prod_{j=2}^{n-m+1}(\Phi_j(q))^{\lfloor n/j\rfloor-\lfloor m/j\rfloor-\lfloor(n-m+1)/j\rfloor}$,
$g(q)=\prod_{j=n-m+2}^n(\Phi_j(q))^{\lfloor n/j\rfloor-\lfloor m/j\rfloor-\lfloor(n-m)/j\rfloor}$. Then $f(q)$, $g(q)$ are monic
 polynomials in $\mathbb{Z}[q]$ and ${\rm deg}(g(q))\geq 1$ because  $\Phi_i(q)$ is a factor of $g(q)$.  By (\ref{h(q)define}), we have $h(q)=f(q)/g(q)$.
Recall that $\Phi_j(q)$, $j=1,\ldots, n$, are irreducible polynomials in $\mathbb{Z}[q]$. We have $g(q)\nmid f(q)$.
 By the polynomial division algorithm,  $f(q)=g(q)f_1(q)+r(q)$, where $f_1(q),r(q)\in \mathbb{Z}[q]$, $r(q)\neq 0$ and ${\rm deg}(r(q))<{\rm deg}(g(q))$.
  Thus, $h(q)=f_1(q)+r(q)/g(q)$. Clearly, if $q$ is a sufficiently large integer (i.e., there is a fixed positive integer $c_{n,m}$ such that $q\geq c_{n,m}$),
 then $h(q)$ is not an integer. Thus, Lemma~\ref{tte54gdg} implies that  $J_q(n,m)$ is a core if $q$ is a sufficiently large integer.  $\qed$

When $m$ and $n-m+1$ are not relatively prime, we guess that $J_q(n,m)$ is a core for all $q$ (which is a power of a prime).

\section*{Acknowledgement}

We are grateful to the referees for  useful comments and suggestions.
This research was   supported by National Natural Science Foundation of China
(Projects 11371072, 11301270, 11271047, 11371204).

\end{document}